\documentclass[a4paper,12pt]{amsart}
\usepackage{amssymb}
\usepackage{ifthen}
\usepackage{graphicx}
\usepackage{amsfonts}
\usepackage{amscd}
\usepackage{amsxtra}
\usepackage{color}

\newtheorem{thm}{Theorem}
\newtheorem{cor}{Corollary}
\newtheorem{lemma}{Lemma}

\newtheorem{conj}{Conjecture}
\theoremstyle{definition}
\newtheorem{example}{Example}
\newtheorem{prob}{Problem}

\newcounter {own}
\def\theown {\thesection       .\arabic{own}}

\newenvironment{rem}{%
\bigskip
\noindent \textsl{{\sl Remark. }}}{\bigskip}
\newenvironment{rems}{%
\bigskip
\noindent \textsl{{\sl Remarks. }}}{\bigskip}

\newenvironment{pf}[1][]{%
 \vskip 3mm
 \noindent
 \ifthenelse{\equal{#1}{}}%
  {{\slshape Proof. }}%
  {{\slshape #1.} }%
 }%
{\qed\bigskip}

\newcounter{alphabet}
\newcounter{tmp}

\makeatletter
\newcommand{\Ref}[1]{\@ifundefined{r@#1}{}{\setcounter{tmp}{\ref{#1}}\Alph{tmp}}}
\makeatother

\newcommand{\IC}{{\mathbb C}}
\newcommand{\ID}{{\mathbb D}}

\def\be{\begin{equation}}
\def\ee{\end{equation}}

\newcommand{\bee}{\begin{enumerate}}
\newcommand{\eee}{\end{enumerate}}

\newcommand{\blem}{\begin{lemma}}
\newcommand{\elem}{\end{lemma}}
\newcommand{\bthm}{\begin{thm}}
\newcommand{\ethm}{\end{thm}}
\newcommand{\bcor}{\begin{cor}}
\newcommand{\ecor}{\end{cor}}
\newcommand{\beg}{\begin{example}}
\newcommand{\eeg}{\end{example}}
\newcommand{\begs}{\begin{examples}}
\newcommand{\eegs}{\end{examples}}
\newcommand{\bdefe}{\begin{defin}}
\newcommand{\edefe}{\end{defin}}
\newcommand{\bprob}{\begin{prob}}
\newcommand{\eprob}{\end{prob}}
\newcommand{\bei}{\begin{itemize}}
\newcommand{\eei}{\end{itemize}}

\newcommand{\bcon}{\begin{conj}}
\newcommand{\econ}{\end{conj}}
\newcommand{\bcons}{\begin{conjs}}
\newcommand{\econs}{\end{conjs}}
\newcommand{\bprop}{\begin{propo}}
\newcommand{\eprop}{\end{propo}}
\newcommand{\br}{\begin{rem}}
\newcommand{\er}{\end{rem}}
\newcommand{\brs}{\begin{rems}}
\newcommand{\ers}{\end{rems}}
\newcommand{\bo}{\begin{obser}}
\newcommand{\eo}{\end{obser}}
\newcommand{\bos}{\begin{obsers}}
\newcommand{\eos}{\end{obsers}}
\newcommand{\bpf}{\begin{pf}}
\newcommand{\epf}{\end{pf}}
\newcommand{\ba}{\begin{array}}
\newcommand{\ea}{\end{array}}
\newcommand{\beq}{\begin{eqnarray}}
\newcommand{\beqq}{\begin{eqnarray*}}
\newcommand{\eeq}{\end{eqnarray}}
\newcommand{\eeqq}{\end{eqnarray*}}

\def\cc{\setcounter{equation}{0}   
\setcounter{figure}{0}\setcounter{table}{0}}

\newcounter{minutes}
\divide\time by 60
\newcounter{hours}
\multiply\time by 60 \addtocounter{minutes}{-\time}

\begin{document}
\bibliographystyle{amsplain}
\title[Where is $f(z)/f'(z)$ univalent?]{Where is $f(z)/f'(z)$ univalent?}

\thanks{
File:~\jobname .tex,
          printed: \number\day-\number\month-\number\year,
          \thehours.\ifnum\theminutes<10{0}\fi\theminutes}

\author[M. Obradovi\'{c}]{Milutin Obradovi\'{c}}
\address{M. Obradovi\'{c},
Department of Mathematics,
Faculty of Civil Engineering, University of Belgrade,
Bulevar Kralja Aleksandra 73, 11000
Belgrade, Serbia. }
\email{obrad@grf.bg.ac.rs}

\author[S. Ponnusamy]{Saminathan Ponnusamy $^\dagger $
}
\address{S. Ponnusamy,
Indian Statistical Institute (ISI), Chennai Centre, SETS (Society
for Electronic Transactions and Security), MGR Knowledge City, CIT
Campus, Taramani, Chennai 600 113, India. }
\email{samy@isichennai.res.in, samy@iitm.ac.in}

 \author[K.-J. Wirths]{Karl-Joachim Wirths}
\address{K.-J. Wirths, Institut f\"ur Analysis und Algebra, TU Braunschweig,
38106 Braunschweig, Germany}
\email{kjwirths@tu-bs.de}

\keywords{Analytic, univalent,  starlike functions, radius of univalence. \\
$
^\dagger$ {\tt  Corresponding author.
}
}

\begin{abstract}
Let ${\mathcal S}$ denote the family of all univalent functions $f$
in the unit disk $\ID$ with the normalization $f(0)=0= f'(0)-1$.
There is an intimate relationship between the operator $P_f(z)=f(z)/f'(z)$
and the Danikas-Ruscheweyh operator $T_f:=\int_{0}^{z}(tf'(t)/f(t))\,dt$.
In this paper we mainly consider the univalence problem of $F=P_f$, where $f$
belongs to some subclasses of ${\mathcal S}$. Among several sharp results and non-sharp results,  
we also show that if $f\in {\mathcal S}$, then  $F \in {\mathcal U}$ in the disk $|z|<r$ with 
$r\leq r_6\approx 0.360794$ and conjecture that the upper bound for such  $r$ is $\sqrt{2}-1$.
\end{abstract}

\maketitle \pagestyle{myheadings}
\markboth{M. Obradovi\'{c}, S. Ponnusamy and K.-J. Wirths}{Where is $f(z)/f'(z)$ univalent?}


\cc
\section{Introduction and Main Results}
Let ${\mathcal B}$ denote the class of analytic functions $\omega (z)$ in
the unit disk $\ID:= \{z\in \IC:\, |z| <1 \}$ such that $\omega(0)=0$ and $|\omega(z)|<1$ for $z\in\ID$.
If $f,g$ are two analytic functions in $\ID$, then we say that $f$ is subordinate to $g$,
written $f\prec g$ or $f(z)\prec g(z),$ if there exists an $\omega \in \mathcal{B}$ such that $f(z)=g(\omega (z)).$
We also note that if $g$ is univalent, then it is easy to show that $f\prec g$ if and only if $f(0)=g(0)$ and $f(\ID)\subset g(\ID).$

We consider the family ${\mathcal A}$ of all functions $f$ analytic in $\ID$
with the normalization $f(0)=0= f'(0)-1$. By ${\mathcal S}$, ${\mathcal S}\subset {\mathcal A}$,
we denote the class of univalent functions in $\ID$.
Certain special subclasses of $\mathcal S$ possess various remarkable features due to their geometrical
properties. By ${\mathcal C}$, ${\mathcal K}$, and ${\mathcal S}^{\star}$
we denote the subclasses of ${\mathcal S}$ which consist of convex, close-to-convex, and
starlike functions, respectively.
For $\beta \in [0,1)$, let $\mathcal{S}^{\star} (\beta)$
denote the usual normalized class of all (univalent) starlike functions of order $\beta $.
Analytically, $f\in \mathcal{S}^{\star} (\beta)$ if $f\in \mathcal{A}$ and satisfies the condition
$$ \frac{zf'(z)}{f(z)} \prec \frac{1+(1-2\beta )z}{1-z}, \quad z\in \ID.
$$
It is well-known that $\mathcal{C}\subsetneq \mathcal{S}^{\star} (1/2)$, and $\mathcal{S}^{\star}:= \mathcal{S}^{\star} (0)$.
At this point it is interesting to note that a function belonging to $\mathcal{S}^{\star} (1/2)$
may not be convex in $|z|<R$ for any $R> \sqrt{2\sqrt{3} -3} = 0.68\ldots $, see \cite[Theorem 1]{Ma2}. We say that $f\in \mathcal{A}$ is
starlike in $|z|<r$ (i.e. to say $f\in \mathcal{S}^{\star}$ in $|z|<r$) for some $0<r\leq 1$, if $f(|z|<r)$ is
starlike with respect to the origin. This means that the last subordination condition is satisfied for $|z|<r$ instead
of the full disk $|z|<1$. Similar convention will be followed for
other classes. We refer to \cite{Du,Go,Pomm} for a detailed discussion on these classes.
Also let us introduce some notations and definitions as follows:
\beqq
{\mathcal U} &=& \left \{ f \in {\mathcal A} :\, \left |U_f(z) \right | < 1
~\mbox{ for $z\in \ID$} \right \},~~ U_f(z)=f'(z)\left (\frac{z}{f(z)} \right )^{2}-1,\\
{\mathcal C}(-1/2) &=& \left \{f\in {\mathcal A}:\,{\rm Re\,}\left(1+\frac{zf''(z)}{f'(z)}\right)>-\frac{1}{2},
~\mbox{ for $z\in \ID$} \right \}, \mbox{ and }\\
{\mathcal G} &=& \left \{f\in {\mathcal A}:\, {\rm Re \,}\left(1+\frac{zf''(z)}{f'(z)}\right)< \frac{3}{2},
 ~\mbox{ for $z\in \ID$}\right \}.
\eeqq
According to Aksent\'ev's theorem \cite{Aks58} (see also \cite{OzNu72}),
the strict inclusion $\mathcal{U}\subsetneq \mathcal{S}$ holds. Moreover, ${\mathcal C}(-1/2) \subset {\mathcal K}$, and
functions in ${\mathcal G}$ are proved to be starlike in $\ID$, see for eg. \cite[Example 1, Equation (16)]{PoRa95}.
See also \cite{ObPoWi2013} for further details and investigation on the class ${\mathcal G}$.

This article concerns with the operator
\be\label{3-13eq1}
F(z):=P_f(z)=\frac{f(z)}{f'(z)}
\ee
for locally univalent functions $f\in\mathcal A$. The main problem is to consider the univalency and starlikeness of $P_f$
when $f$ belongs to some of the subclasses of $\mathcal S$ defined above.

Among others our interest in the operator $P_f$ arose from the fact that
there exists an intimate relation between this one and the Danikas-Ruscheweyh (\cite{DR-1999})
operator
\be\label{6-13eq5}
T_f(z):=\int_{0}^{z}\frac{tf'(t)}{f(t)}\,dt=z+\sum_{n=1}^{\infty}\frac{n}{n+1}c_{n}(f)z^{n+1} \quad (f\in {\mathcal S}),
\ee
where $c_{n}(f)$ $(n\geq 1)$ denote the logarithmic coefficients of $f\in {\mathcal S}$ defined by
$$\log\frac{f(z)}{z}=\sum_{n=1}^{\infty}c_{n}(f)z^{n}.
$$
The conjecture that $T_f\in {\mathcal S}$ for each $f\in {\mathcal S}$ remains open.

The relation between \eqref{3-13eq1} and \eqref{6-13eq5} becomes obvious, when one considers the equivalent operators in the $w$-plane where $w=f(z)$.
Let $g(w)=f^{-1}(w)$ be the function inverse to $f$. If we transform the
operator $P_f$ to the $w$-plane, we get the operator
$$Q(g)(w)=wg'(w)=q(w).
$$
A similar consideration concerning the Danikas-Ruscheweyh operator results in
$$S(g)(w)=\int_0^w \frac{g(u)}{u}\, du=s(w).
$$
Now it is immediately seen that
$$Q^{-1}(q)(w)=\int_0^w \frac{q(u)}{u}\, du = S(q)(w) ~\mbox{ and }~ S^{-1}(s)(w)=ws'(w)=Q(s)(w).
$$

\section{Preliminaries and two examples}

We remark that if $f\in\mathcal{S}$ then $(z/f(z))\neq 0$ in $\ID$ and hence, $f$ can be represented as Taylor's series of the form
\be\label{eq1}
f(z)=\frac{z}{1+\sum_{n=1}^\infty b_nz^n}.
\ee
According to the well-known Area Theorem \cite[Theorem 11 on p.193 of Vol. 2]{Go}, for $f\in \mathcal{S}$ of the form \eqref{eq1}, one has
\be\label{eq7}
\sum_{n=2}^\infty (n-1)|b_n|^2\leq 1
\ee
but this condition is not sufficient for the univalence of $f$. On the other hand, if $f\in \mathcal{A}$ of the form \eqref{eq1} satisfies the condition
\be\label{eq7a}
\sum_{n=2}^{\infty}(n-1)|b_{n}|\leq 1,
\ee
then $f\in\mathcal{U}$. The condition \eqref{eq7a} is also necessary if $b_n\geq 0$ for $n\geq 1$.
The constant $1$ is the best possible in the sense that if
$$\sum_{n=2}^{\infty}(n-1)|b_{n}|= 1+\varepsilon,
$$
for some $\varepsilon>0$, then there exists an $f$ which is not univalent in $\ID$.

Let us continue the discussion with two examples. Consider
$$ f_{1}(z)=\frac{z(1-\frac{z}{2})}{(1-z)^{2}},~\mbox{ and }~f_{2}(z)=z-\frac{z^{2}}{2}.
$$
Then $f_{1}\in {\mathcal C}(-1/2)$ and $f_{2}\in {\mathcal G}$. Define
$$F_{j}(z)= P_{f_j}(z)=\frac{f_{j}(z)}{f'_{j}(z)}, ~\mbox{ for $j=1,2$,}
$$
so that
$$F_{1}(z)=z-\frac{3}{2}z^{2} +\frac{1}{2}z^{3} ~\mbox{ and }~F_{2}(z)=\frac{z(1-\frac{z}{2})}{1-z}.
$$
\begin{enumerate}
 \item We have that
$$F_{1}'(z)=\frac{3}{2}z^{2}-3z+1=\frac{3}{2}(z-r_{+})(z-r_{-}),~r_{\pm} =1\pm \frac{\sqrt{3}}{3}
$$
and therefore $F_{1}'(r_{-})=0$, where $r_{-}=1-\frac{\sqrt{3}}{3}=0.4226497\ldots$.
We claim that ${\rm Re}\,(F_{1}'(z))>0$ for $|z|<r_{-}$. To do this, we observe that
$${\rm Re}\,(F_{1}'(re^{i\theta}))=3r^{2}\cos^{2}\theta -3r\cos\theta+1-\frac{3}{2}r^{2},
$$
then it is easy to show that ${\rm Re}\,(F_{1}'(re^{i\theta}))>0$ for $-1\leq \cos\theta \leq 1$
and $0\leq r <r_{-}$. It means that $F_{1}$ is univalent in the disc $|z|<r_{-}$.

\item It is a simple exercise to see that $F_{2}\in {\mathcal U}$. In fact,
$$\frac{z}{F_{2}(z)}=\frac{1-z}{1-\frac{z}{2}}= 1-\frac{\frac{z}{2}}{1-\frac{z}{2}} = 1-\frac{z}{2}-
\sum_{n=2}^{\infty}b_n z^{n}, \quad b_n=\frac{1}{2^{n}},
$$
so that $z/F_2(z)$ is non-vanishing in $\ID$ and thus,
$$ -z\left(\frac{z}{F_2(z)}\right)'+\frac{z}{F_2(z)} -1 = \left (\frac{z}{F_2(z)}\right )^2F_2'(z)-1=
\left ( \frac{\frac{z}{2}}{1-\frac{z}{2}}\right )^2
$$
from which we easily see that $\left |U_{F_2}(z) \right |<1$
for $z\in \ID$. Indeed, by a direct computation, we see that the function
$w=(z/2)/(1- (z/2))$ maps $\ID$ onto the disk $|w-(1/3)|<2/3$ so that $w\in \ID$ and thus, $w^2\in \ID$.
This observation gives that $\left |U_{F_2}(z) \right |<1$ in $\ID$ and hence, $F_{2}\in {\mathcal U}$. Alternately,
using the series expansion for $F_2$, we find that
$$\sum_{n=2}^{\infty}(n-1)|b_n|= \sum_{n=2}^{\infty}(n-1)\frac{1}{2^{n}}=1
$$
and, by the sufficient condition \eqref{eq7a}, it follows that $F_{2}\in {\mathcal U}$.
\end{enumerate}

%


\section{Main results}
Let $\omega\in {\mathcal B}$.
Then by the Schwarz lemma it follows that
$|\omega(z)|\leq |z|$ for $z\in\ID$ and by the Schwarz-Pick lemma we have
\be\label{6-13eq1}
|\omega'(z)|\leq\frac{1-|\omega(z)|^{2}}{1-|z|^{2}}~\mbox{ for $z\in\ID$}.
\ee
Clearly, $ \frac{\omega(z)}{z}$ is analytic in $\ID$ and $|\omega(z)/z|\leq 1$ in $\ID$.
The Schwarz-Pick lemma, namely, \eqref{6-13eq1}, applied to $\omega(z)/z$ shows that
\be\label{6-13eq2}
|z\omega'(z)-\omega(z)|\leq\frac{|z|^{2}-|\omega(z)|^{2}}{1-|z|^{2}}.
\ee
These three inequalities will be used frequently in the proof of our main results.

\bthm\label{3-13th1}
If $f \in {\mathcal S}^{\star}(\beta )$, then $P_{f}\in {\mathcal U} $ in the disk
 $|z|<1/(1+\sqrt{2(1-\beta)})$. The result is sharp (as for univalence) as the function $z/(1-z)^{2(1-\beta)}$ shows.
\ethm\bpf
Each $f \in {\mathcal S}^{\star}(\beta )$ and $F=P_f$ defined by \eqref{3-13eq1} can be written as
$$\frac{zf'(z)}{f(z)}=\frac{1+(1-2\beta)\omega(z)}{1-\omega(z)} ~\mbox{ and }~
F(z)=\frac{z(1-\omega(z))}{1+(1-2\beta)\omega(z)},
$$
where $\omega \in {\mathcal B}$.
Clearly, $ \frac{\omega(z)}{z}$ is analytic in $\ID$ and $|\omega(z)/z|\leq 1$ in $\ID$.
Using the last two relations, we observe that
\be\label{6-13eq3}
U_{F}(z) = -z\left(\frac{z}{F(z)}\right)'+\frac{z}{F(z)} -1
=\frac{zf'(z)}{f(z)}-z\left(\frac{zf'(z)}{f(z)} \right)'-1
\ee
and thus,
$$U_{F}(z) = 2(1-\beta)\left (\frac{\omega(z)}{1-\omega(z)}-\frac{z\omega'(z)}{(1-\omega(z))^{2}}\right)
= 2(1-\beta)\left (\frac{(\omega(z)-z\omega'(z))-\omega^2(z)}{(1-\omega(z))^{2}}\right) 
$$
from which and \eqref{6-13eq2}, we obtain that
\beqq
|U_{F}(z)|
&\leq& 2(1-\beta) \left (\frac{|\omega(z)-z\omega'(z)|}{(1-|\omega(z)|)^{2}} + \frac{|\omega(z)|^{2}}{(1-|\omega(z)|)^{2}}\right )\\
&\leq & 2(1-\beta) \left (\frac{\frac{|z|^{2}-|\omega(z)|^{2}}{1-|z|^{2}}}{(1-|\omega(z)|)^{2}} + \frac{|\omega(z)|^{2}}{(1-|\omega(z)|)^{2}}\right )\\
&=&\frac{2(1-\beta)|z|^{2}}{1-|z|^{2}}\left (\frac{1+|\omega(z)|}{1-|\omega(z)|}\right )
\leq \frac{2(1-\beta)|z|^{2}}{1-|z|^{2}}\left (\frac{1+|z|}{1-|z|}\right ) =\frac{2(1-\beta)|z|^{2}}{(1-|z|)^{2}}
\eeqq
which can easily seen to be less than $1$ if $|z|<1/(1+\sqrt{2(1-\beta)}).$ Thus, $F$ belongs to $\mathcal U$ in the disk $|z|<1/(1+\sqrt{2(1-\beta)})$.

To prove the sharpness part, we consider $k_{\beta}(z)=z/{(1-z)^{2(1-\beta)}}$
and define
$$F_{\beta}(z)= P_{k_{\beta}}(z)=\frac{k_{\beta}(z)}{k_{\beta}'(z)}.
$$
Then we see that $k_{\beta}\in {\mathcal S}^*(\beta)$ and
$$F_{\beta}(z) =\frac{z(1-z)}{1+(1-2\beta)z}
~\mbox{ and }~\frac{z}{F_{\beta}(z)}=\frac{1+(1-2\beta)z}{1-z}=1+2(1-\beta)\sum_{n=1}^{\infty}z^{n}.
$$
Define $G_\beta (z)=\frac{1}{r}F_{\beta}(r z)$ and observe that
$$\frac{z}{G_{\beta}(z) }=1+2(1-\beta)\sum_{n=1}^{\infty}r^{n}z^{n}.
$$
According to \eqref{eq7a},
the function $G_{\beta}$ is in $\mathcal U$ (and hence is univalent in $\ID$) if and only if
$$2(1-\beta)\sum_{n=2}^{\infty}(n-1)r^{n}\leq 1, ~\mbox{ i.e. }~ \frac{2(1-\beta)r^{2}}{(1-r)^{2}}\leq 1.
$$
The gives the condition $0<r\leq r_{1}=1/(1+\sqrt{2(1-\beta)})$. Thus, the function
$F_{\beta}$ is univalent in the disk $|z|<r_1$ and not in any larger larger disk with center at the origin. Note also that
$$F_\beta'(z)=\frac{1-2z-(1-2\beta)z^2}{(1+(1-2\beta)z)^2}
$$
and thus, $F_\beta '(r_1)=0$. Moreover,
$$U_{F_\beta}(z)=\frac{1-2z-(1-2\beta)z^2}{(1-z)^2} -1
$$
showing that $U_{F_\beta}(r_1)=-1$. Thus, the number $r_1$ is best both for univalence and also for $\mathcal U$.
The proof is complete.
\epf

\bcor\label{3-13cor1}
If $f \in {\mathcal S}^{\star}$, then $P_{f}\in {\mathcal U}\cap {\mathcal S}^{\star}$ in the disk
$|z|<\sqrt{2}-1$. The result is sharp (as for univalence) as the Koebe function $z/(1-z)^{2}$ shows.
\ecor
\bpf
It suffices to prove the starlikeness part since $P_{f}\in {\mathcal U}$ follows from Theorem \ref{3-13th1} by taking $\beta =0$.
Thus, for the proof of the second part, it suffices to observe by \eqref{6-13eq1} that
$$\left|\frac{zF'(z)}{F(z)}-1\right|=\left|-\frac{2z\omega'(z)}{1-\omega^{2}(z)}\right|
\leq \frac{2|z|\,|\omega'(z)|}{1-|\omega(z)|^{2}}\leq \frac{2|z|}{1-|z|^{2}}
$$
which is again less than $1$ provided $|z|<\sqrt{2}-1$. In particular, $F$ is starlike in the disk $|z|<\sqrt{2}-1$.
Sharpness part follows from the discussion in Theorem \ref{3-13th1} with $\beta =0$.
\epf

\bcor\label{3-13th2}
If $f \in {\mathcal S}^{\star}(1/2)$, then $P_{f}\in {\mathcal U}\cap {\mathcal S}^{\star}$ in the disk $|z|<1/2$. The result is sharp
as the function $z/(1-z)$ shows.
\ecor\bpf
Choose $\beta =1/2$ in Theorem \ref{3-13th1} and observe that it suffices to prove the starlikeness part.
As in the proof of Theorem \ref{3-13th1}, for each $f \in {\mathcal S}^{\star}(1/2)$, we have
$$\frac{zf'(z)}{f(z)}=\frac{1}{1-\omega(z)} ~\mbox{ and }~ F(z)= z(1-\omega(z))
$$
for some $\omega \in {\mathcal B}$. 
%
By \eqref{6-13eq1} and the fact that $|\omega(z)|\leq |z|$, we obtain
$$\left|\frac{zF'(z)}{F(z)}-1\right|= \left|\frac{-z\omega'(z)}{1-\omega(z)}\right|
\leq \frac{|z|\,|\omega'(z)|}{1-|\omega(z)|} \leq\frac{|z|(1+|\omega(z)|)}{1-|z|^{2}}
\leq \frac{|z|(1+|z|)}{1-|z|^{2}}=\frac{|z|}{1-|z|}
$$
which is less than $1$ if $|z|<1/2$. Note that for $f(z)=z/(1-z)$, one has $F(z)= z-z^{2}$ and thus, $|F'(z)-1|=2|z|<1$ for $|z|<1/2$ and
$F'(1/2)=0$. Thus, $F$ is univalent in the disk $|z|<1/2$ and not in any larger disk with center at the origin. Also, it is easy to see that $F(z)$ is starlike
for $|z|<1/2$. The desired conclusion follows.
\epf

\bcor\label{3-13cor3a}
If $f \in {\mathcal S}^{\star}(1/2)$ such that $f''(0)=0$, then $P_f$ is starlike in the disk $|z|<r_2$, where
$r_2\approx 0.543689$ is the root of the equation $\phi_2 (r)=0$, where
$$ \phi_2 (r)=r^{3}+r^{2}+r-1.
$$
\ecor\bpf
Clearly, we just need to apply Corollary \ref{3-13th2} with $|\omega(z)|\leq|z|^{2}$. This will lead to the
inequality
$$ \left|\frac{zF'(z)}{F(z)}-1\right|\leq \frac{|z|(1+|z|^{2})}{1-|z|^{2}}
$$
which is clearly less than $1$ if $|z|^{3}+|z|^{2}+|z|-1<0$. The result follows.
\epf

\bcor
Let $f$ belong to either ${\mathcal S}^{\star}(1/2)$ or ${\mathcal C}(-1/2)$, such that
$f''(0)=0$. Then $F\in {\mathcal U}$ in the disk $|z|<1/\sqrt{3}$.
\ecor
\bpf
It known that \cite[p.~68]{MiMo2} if ${\mathcal C}(-1/2)$ with $f''(0)=0$, then $f\in {\mathcal S}^{\star}(1/2)$. In view of this result,
it suffices to prove the corollary when $f$ belongs to ${\mathcal S}^{\star}(1/2)$ with $f''(0)=0$. However, using the proof of
Theorem \ref{3-13th1} with $\beta =1/2$ and $|\omega(z)|\leq|z|^{2}$, we easily obtain that
$$|U_{F}(z)| \leq \frac{|z|^{2}}{1-|z|^{2}}\left (\frac{1+|\omega(z)|}{1-|\omega(z)|}\right )\leq
\frac{|z|^{2}}{(1-|z|^{2})}\left (\frac{1+|z|^2}{1-|z|^2} \right )
$$
which is less than $1$ provided $1-3|z|^2>0$ and this gives the disk $|z|<1/\sqrt{3}$. The proof is complete.
\epf

A locally univalent function $f\in {\mathcal A}$ is said to
belong to ${\mathcal G}(\alpha )$, for some $\alpha \in (0,1]$, if it satisfies the condition
\begin{equation}\label{Ob-S3-11-eq6a}
{\rm Re} \left (1+\frac{zf''(z)}{f'(z)}\right )<1+\frac{\alpha }{2},
\quad \mbox{$z\in {\mathbb D}$}.
\end{equation}
Thus, we have ${\mathcal G}:={\mathcal G}(1)$.

\bthm \label{3-13th3}
If $f \in {\mathcal G}(\alpha )$ for some $\alpha \in (0,1]$, then $P_f$ is starlike in the disk $|z|<1+\alpha -\sqrt{\alpha (1+\alpha )}$ .
\ethm
\bpf Let $f \in {\mathcal G}(\alpha )$ and $F$ be given by \eqref{3-13eq1}. Then we have (see eg. \cite[Theorem 1]{JoOb95})
$$ \frac{zf'(z)}{f(z)} \prec \frac{(1+\alpha )(1-z)}{1+\alpha -z}, \quad \mbox{$z\in {\mathbb D}$,}
$$
and thus,
we may write
$$\frac{zf'(z)}{f(z)}=\frac{(1+\alpha )(1-\omega(z))}{1+\alpha -\omega(z)}
 ~\mbox{ and }~ F(z)= P_f= \frac{z(1+\alpha -\omega(z))}{(1+\alpha ) (1-\omega(z))}
$$
for some $\omega \in {\mathcal B}$. By a computation, we obtain that
$$\frac{zF'(z)}{F(z)}-1 =\frac{\alpha z\omega'(z)}{(1-\omega(z))(1+\alpha -\omega(z))}
$$
and, as before, it follows from the Schwarz-Pick lemma that
$$\left|\frac{zF'(z)}{F(z)}-1\right| \leq \frac{\alpha |z|\,|\omega'(z)|}{(1+\alpha -|\omega(z)|)(1-|\omega(z)|)}
\leq \frac{\alpha |z|}{(1+\alpha -|z|)(1-|z|)}
$$
which is less than $1$ provided $\phi _3 (|z|)>0$, where $ \phi _3 (r) =r^2 -2(1+\alpha )r +1+\alpha .$
Thus, we conclude that $P_f$ is starlike in the disk $|z|<r_3(\alpha)=1+\alpha -\sqrt{\alpha (1+\alpha )}$, where
$r_3(\alpha)$ is the root of the equation $\phi _3 (r)=0$ in the interval $(0,1]$. The theorem follows.
\epf


Taking $\alpha =1$ gives

\bcor \label{3-13cor3}
If $f \in {\mathcal G}$, then $P_f$ is starlike in the disk $|z|<2-\sqrt{2} \approx 0.585786$ .
\ecor


The same reasoning gives as in Corollary \ref{3-13cor3a} the following.
\bcor
If $f \in {\mathcal G}(\alpha)$ such that $f''(0)=0$ and for some $\alpha \in (0,1]$, then $P_f$ is starlike in $|z|<r_4(\alpha)$, where
$r_4(\alpha)$ is the root in the interval $(0,1]$ of the equation $\phi_4 (r)=0$,
$$\phi_4 (r)= r^{4}-\alpha r^{3}-(2+\alpha)r^{2}-\alpha r +1+\alpha .
$$
\ecor
\bpf
In this case, the corresponding inequality for $f\in {\mathcal G}(\alpha)$ in Theorem \ref{3-13th3} becomes
$$\left|\frac{zF'(z)}{F(z)}-1\right|
\leq \frac{\alpha |z|}{1-|z|^{2}}\left (\frac{1+ |\omega(z)|}{1+\alpha -|\omega(z)|}\right )
\leq \frac{\alpha |z|}{1-|z|^{2}}\left (\frac{1+|z|^{2}}{1+\alpha -|z|^{2}}\right )
$$
which is less than $1$ if $\phi_4 (|z|)>0$. The result follows.
\epf

Setting $\alpha=1$ gives

\bcor
If $f \in {\mathcal G}$ such that $f''(0)=0$, then $P_f$ is starlike in $|z|<r_4$, where
$r_4\approx 0.64731$ is the root in the interval $(0,1]$ of the equation $r^{4}-r^{3}-3r^{2}-r+2=0 $.
\ecor

\bthm\label{3-13th5}
If $f \in {\mathcal G}(\alpha)$ for some $\alpha \in (0,1]$, then $F\in {\mathcal U}$
in the disk $|z|<r_5(\alpha )$, where $ r_5(\alpha )=\sqrt{\frac{-\alpha +\sqrt{(1+\alpha)^2+1}}{2}}$.
\ethm
\bpf
Let $f \in {\mathcal G}(\alpha )$ and $F=P_f$ be given by \eqref{3-13eq1}. Then, following the proof of Theorem \ref{3-13th3}, one has
$$ \frac{z}{F(z)}-1 = -\frac{\alpha \omega(z)}{1+\alpha -\omega(z)}
$$
and, using this relation, we find that
\beqq
U_{F}(z) &=& -\frac{\alpha \omega(z)}{1+\alpha -\omega(z)} +\frac{\alpha (1+\alpha)z\omega'(z)}{(1+\alpha -\omega(z))^{2}}\\
&=&\frac{\alpha [(1+\alpha)(z\omega'(z) -\omega (z)) + \omega ^2(z)]}{(1+\alpha -\omega(z))^{2}}
\eeqq
so that, by \eqref{6-13eq2}, we easily have as before that
\beqq
|U_{F}(z)|
&\leq& \frac{\alpha}{(1+\alpha -|\omega(z)|)^{2}} \left ( (1+\alpha)\left (\frac{|z|^{2}-|\omega(z)|^{2}}{1-|z|^{2}} \right ) +|\omega(z)|^{2}\right )\\
&=& \frac{\alpha}{1-|z|^{2}}\left (\frac{-(\alpha +|z|^{2})|\omega(z)|^{2}+(1+\alpha)|z|^{2}}{(1+\alpha-|\omega(z)|)^{2}}\right )
=\frac{\alpha \phi(t)}{1-r^{2}} ,
\eeqq
where we put $|z|=r$, $|\omega(z)|=t$ and
$$\phi(t)=\frac{-(\alpha +r^{2})t^{2}+(1+\alpha )r^{2}}{(1+\alpha -t)^{2}}, \,\,0\leq t\leq r.
$$
We compute that
$$\phi'(t)=\frac{2(1+\alpha )}{(1+\alpha -t)^{3}}\left[-(\alpha +r^{2})t+r^{2}\right],
$$
and it is easy to see that $\phi$ attains its maximum value $\phi (t_0)$, where $t_{0}=\frac{r^{2}}{\alpha+r^{2}}$
and $\phi '' (t_0)<0$.
A calculation gives
$$\phi(t_{0})=\frac{r^{2}(\alpha +r^{2})}{\alpha (1+\alpha +r^{2})}
$$
and thus, we have
$$|U_{F}(z)|\leq \frac{\alpha \phi(t_{0}) }{1-r^{2}} =\frac{r^{2}(\alpha +r^{2})}{(1-r^{2})(1+\alpha+r^{2})}
$$
which is less than $1$ if $2r^{4}+2\alpha r^{2}-(1+\alpha)<0$. This gives that $|U_{F}(z)|<1$ for $0<r\leq r_5(\alpha )$,
where $r_5(\alpha )$ is the root of the equation $2r^{4}+2\alpha r^{2}-(1+\alpha)=0$, that lies in the interval $(0,1)$. The conclusion follows.
\epf

The choice $\alpha =1$ yields the following.

\bcor
If $f \in {\mathcal G}$, then $F$ belongs to the class ${\mathcal U}$
in the disk $|z|<\sqrt{\frac{\sqrt{5}-1}{2}}\approx 0.78615$.
\ecor

\bthm\label{3-13th6}
Let $f\in {\mathcal S}$ with $a_2=f''(0)/2!$. Then $F$ belongs to ${\mathcal U}$ in the disk $|z|<r_6(|a_{2}|)$, where $r_6(|a_{2}|)$
is the root of the equation $\phi _5(r)=0$ that lies in the interval $(0,1)$, where
$$\phi _5(r)=(a+1-\frac{1}{4}b^2)r^{10}-(5a+5-\frac{5}{4}b^2)r^8+
(19a+10-\frac{19}{4}b^2)r^6
+(9a-10-\frac{9}{4}b^2)r^4
+5r^{2}-1
$$
with $b=|a_2|$ and $a=\frac{2\pi^{2}-12}{3} \approx 2.57974.$
\ethm
\bpf
Let $f\in {\mathcal S}$ and following the idea of \cite[Theorem 4]{obpo-2014}, we consider
\be\label{6-13eq4}
\log\frac{f(z)}{z}=\sum_{n=1}^{\infty}c_{n}(f)z^{n},
\ee
where $c_{n}(f)$ $(n\geq 1)$ denote the logarithmic coefficients of $f$
with $c_1(f)=a_2$. Further, for $f\in {\mathcal S}$ the following sharp
inequality is known from the work of Roth \cite[Theorem 1.1]{Roth-07}
$$\sum_{n=1}^{\infty}\left(\frac{n}{n+1}\right)^{2}|c_{n}(f)|^{2}\leq \frac{2\pi^{2}-12}{3}
=a .
$$
By \eqref{6-13eq4}, we obtain
$$\frac{zf'(z)}{f(z)}-1 =\sum_{n=1}^{\infty}nc_{n}(f)z^{n}
$$
which by the relation \eqref{6-13eq3} gives that
$$U_{F}(z)=-\sum_{n=1}^{\infty}n(n-1)c_{n}(f)z^{n}
$$
and thus, by the Cauchy-Schwarz inequality, we obtain that
\beqq
|U_{F}(z)|& = & 
\left|\sum_{n=2}^{\infty} n(n-1)c_{n}(f)z^{n} \right|\\
& \leq & \left( \sum_{n=2}^{\infty}\left(\frac{n}{n+1}\right)^{2}|c_{n}(f)|^{2}\right)^{\frac{1}{2}}
\left( \sum_{n=2}^{\infty}(n^{2}-1)^{2}|z|^{2n}\right)^{\frac{1}{2}}\\
& \leq & \left( a-\frac{1}{4}|c_{1}(f)|^{2}\right)^{\frac{1}{2}}
\left(\frac{|z|^{4}(|z|^{6}-5|z|^{4}+19|z|^{2}+9)}{(1-|z|^{2})^{5}}\right)^{\frac{1}{2}}
\eeqq
which is less than $1$ whenever,
$$\left( a-\frac{1}{4}|c_{1}(f)|^{2}\right)|z|^{4}(|z|^{6}-5|z|^{4}+19|z|^{2}+9) <(1-|z|^{2})^{5}.
$$
If we put $r=|z|$, then the last inequality is equivalent to $\phi _5(r):=\phi _5(r, |a_2|)<0$, where $\phi _5(r)$ is as in the statement.
The desired result follows.
\epf

\bcor\label{3-13cor4}
Let $f\in {\mathcal S}$ with $f''(0)=0$, and $a=\frac{2\pi^{2}-12}{3}.$ Then $F$ belongs to ${\mathcal U}$ in the disk $|z|<r_6$,
where $r_6\approx 0.360794$ is the root of the equation
$$(a+1)r^{10}-5(a+1)r^{8}+(19a+10)r^{6}+(9a-10)r^{4}+5r^{2}-1=0,
$$
that lies in the interval $(0,1)$.
\ecor
\bpf Set $a_2=0$ in Theorem \ref{3-13th6}.
\epf

It is a simple exercise to see that the values $r_6(|a_{2}|)$, as the roots of the equation $\phi _5(r)=0$, increase with increasing values of $|a_2|\in [0,2]$.
For a ready reference, we included in Table \ref{table1} a list of values of $r_6(|a_{2}|)$ for certain choices of $|a_2|$.
This observation shows that if $f\in {\mathcal S}$, then $F \in {\mathcal U}$ in the disk $|z|<r$ and the lower bound for $r$
by Corollary \ref{3-13cor4} is $r_6\approx 0.360794$. We end the discussion with a conjecture that the upper bound for the value of $r$ is $\sqrt{2}-1$ which is attained by the Koebe function.

\begin{table}
\begin{center}
\begin{tabular}{|c|c||c|c|}
\hline
{\bf Values of $|a_2|$}&{\bf values of $r_6(|a_{2}|)$}&{\bf Values of $|a_2|$}&{\bf values of $r_6(|a_{2}|)$}\\
\hline
0.25&0.361166& 1.25&0.370874\\
\hline
0.5&0.362294&1.5&0.375923\\
\hline
0.75&0.364226&1.75&0.382504\\
\hline
1&0.367042&2&0.391124\\
\hline
\end{tabular}
\vspace{.2cm}
\caption{Values of $r_6(|a_{2}|)$ for different values of $|a_2|$ \label{table1}}
\end{center}
\end{table}

\subsection*{Acknowledgement}
The work of the first author was supported by MNZZS Grant, No. ON174017, Serbia.
The second author is currently on leave from Indian Institute of Technology Madras, India.


\begin{thebibliography}{150}
\bibitem{Aks58} {\sc L.~A. Aksent\'ev},
\textrm{Sufficient conditions for univalence of regular functions (Russian)},
\textit{Izv. Vys\v s. U\v cebn. Zaved. Matematika}
\textbf{1958}(4)(1958), 3--7.


\bibitem {DR-1999} {\sc N. Danikas and St. Ruscheweyh},
Semi-convex hulls of analytic functions in the unit disk,
\textit{Analysis} {\bf 4}(1999), 309--318.

\bibitem{Du} {\sc P.~L.~Duren},
\emph{Univalent functions} (Grundlehren der
mathematischen Wissenschaften 259, New York, Berlin, Heidelberg, Tokyo), Springer-Verlag, 1983.


\bibitem {Go} {\sc A.~W. Goodman}, \emph{Univalent functions}, Vols. 1-2,
Mariner, Tampa, Florida, 1983.


\bibitem{JoOb95} {\sc I. Jovanovi\'c, and M. Obradovi\'c,}
\textrm{A note on certain classes of univalent functions},
\textit{Filomat}, {\bf 9}(1)(1995), 69--72.



%

\bibitem{obpo-2014} {\sc M. Obradovi\'c and S. Ponnusamy},
Univalence of quotient of analytic functions,
\textit{Appl. Math. Comput.} \textbf{247}(2014), 689--694.

\bibitem{ObPoWi2013} {\sc M. Obradovi\'{c}, S. Ponnusamy and K.-J. Wirths,}
\textrm{Coefficient characterizations and sections for some univalent functions,}
\textit{Sib. Math. J.}  \textbf{54}(1)(2013), 679--696.

%
\bibitem{Ma2} {\sc T.~H.~MacGregor,}
The radius of convexity for starlike functions of order $1/2$,
\textit{Proc. Amer. Math. Soc.} {\bf 14}(1963), 71--76.


\bibitem{MiMo2} {\sc S. S. Miller and P. T. Mocanu}, {\em Differential
subordinations: Theory and Applications}, Marcel Dekker, Inc. New York. Basel, No. {\bf 225}, pp.459 (2000).

\bibitem{OzNu72} {\sc S. Ozaki and M. Nunokawa},
\textrm{The Schwarzian derivative and univalent functions},
\textit{Proc. Amer. Math. Soc.} {\bf 33}(1972), 392--394.

\bibitem{Pomm} {\sc Ch. Pommerenke},
\emph{Univalent functions},
Vandenhoeck and Ruprecht, G\"ottingen, 1975.

\bibitem{PoRa95} {\sc S. Ponnusamy and S. Rajasekaran},
New sufficient conditions for starlike and univalent functions,
\textit{Soochow J. Math.} {\bf 21}(2)(1995), 193--201.

\bibitem{Roth-07} {\sc O. Roth}, A sharp inequality for the logarithmic coefficients of univalent functions,
\textit{Proc. Amer. Math. Soc. } {\bf 138}(7)(2007), 2051--2054.

\end{thebibliography}
\end{document}